\documentclass[12pt]{amsart}
\usepackage{amsfonts}

\theoremstyle{definition}

\theoremstyle{remark}

\newcommand{\const}{\mathop{\rm const}\limits}

\newcommand{\supp}{\mathop{\rm supp}\limits}

\begin{document}

\begin{center}

{\bf CESARO-HARDY OPERATORS ON BILATERAL \\
\vspace{3mm}
 GRAND LEBESGUE SPACES} \par

\vspace{3mm}
{\bf E. Ostrovsky}\\

e-mail: galo@list.ru \\

\vspace{3mm}

{\bf L. Sirota}\\

e-mail: sirota3@bezeqint.net \\

\vspace{3mm}

Department of Mathematics and Statistics, Bar-Ilan University, 59200, Ramat Gan, Israel. \\

\vspace{3mm}

\begin{center}
 Abstract. \\

\end{center}

{\it We obtain in this short article the non-asymptotic estimations for  the norm of (generalized)
Cesaro-Hardy integral operators in the so-called
Bilateral Grand Lebesgue Spaces. We also give examples to show the
sharpness of these inequalities.} \\

\end{center}

\vspace{3mm}


\vspace{3mm}


\vspace{3mm}

\section{Introduction}

\vspace{3mm}

 The linear integral operator $ U_{\alpha, \beta, \lambda} [f](x) = U[f](x), $ or, wore precisely, the
 {\it family of operators} of a view

 $$
U(x) = U_{\alpha,\beta,\lambda}[f](x) = x^{- \beta }\int_{0}^x \frac{ y^{-\alpha} \ f(y) \ dy}{|x - y|^{\lambda}} \eqno(1.0)
 $$
is called {\it generalized Cesaro-Hardy integral operator,} or {\it fractional integral.} \par

 Here  $ x,y \in (0, \infty), \  \alpha,\beta, \lambda = \const \in (0,1),  \ \alpha + \beta + \lambda < 1. $\par

 We denote as usually the classical $ L_p $ Lebesgue - Riesz norm
 $$
 |f|_p = \left( \int_{X} |f(x)|^p \ dx  \right)^{1/p};  \  f \in L_p \
 \Leftrightarrow |f|_p < \infty, \eqno(1.1)
 $$
and  denote $ L(a,b) = \cap_{p \in (a,b)} L_p.  $ \par
 Here $ X = R_+  $ or $ X = (R_+)^d  $ or $ X = R^d.  $  \par
 The case of operators of a view

 $$
I(x) = I_{\alpha,\beta,\lambda}[f](x) = ||x||^{- \beta }\int_{R^d} \frac{ ||y||^{-\alpha} \ f(y) \ dy}{||x - y||^{\lambda}}
 $$
was considered in \cite{Okikiolu1}, chapter 11, see also \cite{Ostrovsky3}; it was proved in particular the following estimation.
 Define

 $$
 p_-^{(d)} = \frac{d}{d-\alpha}, \hspace{4mm} p_+^{(d)} = \frac{d}{d - \alpha - \lambda},
 $$

 $$
 p_- = p_-^{(1)} =  \frac{1}{1-\alpha}, \hspace{4mm} p_+ = p_+^{(1)} = \frac{1}{1 - \alpha - \lambda},
  $$

 $$
 q_-^{(d)}= \frac{d}{\beta + \lambda}, \hspace{4mm} q_+^{(d)} = \frac{d}{\beta},
 $$

 $$
  q_- = q_-^{(1)} = \frac{1}{\beta + \lambda}, \hspace{4mm}  q_+ = q_+^{(1)} = \frac{1}{\beta} ,
 $$

 $$
 \kappa^{(d)} = (\alpha + \beta + \lambda)/d; \ \kappa = \kappa_1 = \alpha + \beta + \lambda.  \eqno(1.2)
 $$

Define also for the arbitrary value $ p$ from the  set $ p \in (p_-^{(d)}, p_+^{(d)})   $ the  correspondent value
$  q  = q(p), \ q \in (q_-^{(d)}, q_+^{(d)})  $ as follows:

$$
1 + \frac{1}{q} = \frac{1}{p} + \frac{\alpha + \beta + \lambda}{d} =  \frac{1}{p} +  \kappa^{(d)}. \eqno(1.3)
$$
 The identity (1.3) defined uniquely  the function $  p = p(q) $ and  inversely the function $ q = q(p). $ \par

 It is proved in \cite{Ostrovsky4} that

 $$
 |I_{\alpha,\beta,\lambda}[f]|_{q(p)} \le V(\alpha, \beta, \lambda; p) \cdot  |f|_p, \  p \in (p_-^{(d)}, p_+^{(d)}), \eqno(1.4)
 $$
 where for the {\it optimal, } i.e. minimal value $ V(\alpha, \beta, \lambda; p): $

$$
V(\alpha, \beta, \lambda; p) \stackrel{def}{=} \sup_{f \in L( p), f \ne 0 }
\left[ \frac{|I_{\alpha,\beta,\lambda}[f]|_{q(p)}}{|f|_p} \right]
$$
are true the following estimates:

 $$
\frac{C_1(\alpha, \beta, \lambda)}{ \left[(p-p_-^{(d)}) \ (p_+^{(d)} - p) \right]^{\kappa^{(d)} }} \le  V(\alpha, \beta, \lambda; p) \le
\frac{C_2(\alpha, \beta, \lambda)}{ \left[(p-p_-^{(d)}) \ (p_+^{(d)} - p) \right]^{\kappa^{(d)} }}, \eqno(1.5)
 $$

\vspace{4mm}

 $$
  C_1(\alpha, \beta, \lambda), \  C_2(\alpha, \beta, \lambda)  \in (0,\infty). \eqno(1.5a)
 $$

\vspace{4mm}

{\bf  Our purpose is the extension of inequality (1.5) into the generalized  operator of Cesaro - Hardy view: }

\vspace{4mm}

 $$
 |U_{\alpha,\beta,\lambda}[f] |_{q(p)} \le K(\alpha, \beta, \lambda; p) \cdot  |f|_p, \  p \in (p_-^{(d)}, p_+^{(d)}), \eqno(1.6)
 $$

{\it with exact "constant" $ K(\alpha, \beta, \lambda; p) $ estimation, alike (1.5) - (1.5a). }\par

 As before, we will understood in the capacity of the coefficient $ K = K(\alpha, \beta, \lambda; p)  $
 its minimal value:

 $$
 K(\alpha, \beta, \lambda; p) \stackrel{def}{=}
 \sup_{f \in L(p),  f \ne 0}  \left[ \frac{|U_{\alpha,\beta,\lambda}[f] |_{q(p)}}{|f|_p}  \right].  \eqno(1.7)
 $$

 Notice that the case $ \alpha + \beta + \lambda = 1  $ was investigated in the classical book
\cite{Hardy1}, with exact constant computation.  Therefore, we do not impose this condition. \par

 We will obtain also the generalization  of these estimations on the so-called Grand Lebesgue Spaces (GLS).
 Note that the Sobolev's weight space  estimates for these operators are obtained in a recent  article
\cite{Lizama1},  without constants  evaluating. \par

\vspace{5mm}

  These operators are used in the theory of Fourier transform, probability theory, theory of PDE,
in the functional analysis, in particular, in the theory of interpolation of operators etc.,
 see for instance \cite{Bennet1},\cite{Carro1}, \cite{Krein1}, \cite{Perez1}, \cite{Stein1}. \par

 One of absolutely unexpected application of these estimations are  in the theory of  Navier-Stokes equations, see e.g.
\cite{Fujita1}, \cite{Giga4},  \cite{Kato1}, \cite{Ostrovsky5}. Authors hope to use further the results of this report
in the  theory of Navier-Stokes equation. \par

\vspace{5mm}

 We use symbols $C(X,Y),$ $C(p,q;\psi),$ etc., to denote positive
constants along with parameters they depend on, or at least
dependence on which is essential in our study. To distinguish
between two different constants depending on the same parameters
we will additionally enumerate them, like $C_1(X,Y)$ and
$C_2(X,Y).$ The relation $ g(\cdot) \asymp h(\cdot), \ p \in (A,B), $
where $ g = g(p), \ h = h(p), \ g,h: (A,B) \to R_+, $
denotes as usually

$$
0< \inf_{p\in (A,B)} h(p)/g(p) \le \sup_{p \in(A,B)}h(p)/g(p)<\infty.
$$
The symbol $ \sim $ will denote usual equivalence in the limit
sense.\par
We will denote as ordinary the indicator function
$$
I(x \in A) = 1, x \in A, \ I(x \in A) = 0, x \notin A;
$$
here $ A $ is a measurable set.\par
 All the passing to the limit in this article may be grounded by means
 of Lebesgue dominated convergence theorem.\par

\bigskip

\section{Main result: upper and lower  estimations for Cesaro-Hardy operator}

\vspace{3mm}

  {\it  We consider in this section only  the one - dimensional case for these operators:  } $  d = 1.$ \par

\vspace{3mm}

{\bf Lemma 2.1.}
 {\it If the inequality (1.6) there holds for every function $  f  $ from the Schwartz space $ S(R_+),  $  then }

$$
1 + \frac{1}{q} = \frac{1}{p} + \alpha + \beta + \lambda =  \frac{1}{p} +  \kappa. \eqno(2.1)
$$

{\bf Proof.} We will use the well-known {\it scaling, or equally, dilation method, } see \cite{Stein1}, \cite{Talenti1}.
Indeed, let the inequality (1.6) be satisfied for {\it some}  function $ f(\cdot) \ne 0 $ from the set $ S(R_+).  $  \par
 Let $  \gamma = \const \in (0,\infty); $ consider the dilation function

 $$
 f_{\gamma}(x) = T_{\gamma}[f](x) = f(\gamma \ x).
 $$
 Evidently, $ f_{\gamma}(\cdot) \in  S(R_+). $ Therefore

 $$
 | U_{\alpha,\beta, \lambda} [ T_{\gamma}[f]]|_q \le K(\alpha,\beta,\lambda;p) \cdot |T_{\gamma} [f]|_p. \eqno(2.2)
 $$

 We get consequently after simple calculations:

 $$
  |T_{\gamma} [f]|_p = \gamma^{-1/p} \ |f|_p,
 $$

 $$
 | U_{\alpha,\beta, \lambda} [ T_{\gamma}[f]]|_q = \gamma^{\alpha + \beta + \lambda -1 - 1/q  } \  | U_{\alpha,\beta, \lambda} [f]|_q.
 $$
  We conclude substituting into (1.6)

$$
\gamma^{\alpha + \beta + \lambda -1 - 1/q  } \cdot  | U_{\alpha,\beta, \lambda} [f]|_q \le K(\alpha,\beta,\lambda;p) \cdot
\gamma^{-1/p} \cdot |f|_p.
$$
 Since $ \gamma  $ is arbitrary positive number, we see

$$
\alpha + \beta + \lambda -1 - 1/q  = -1/p,
$$
which is equivalent to the assertion  of Lemma 1. \par

\vspace{3mm}

{\bf Lemma 2.2.}
 {\it If the inequality (1.6) there holds  with finite value of $ K(\alpha, \beta, \lambda; p)   $
 for every function $  f  $ from the space $ L_p(R_+),  $  then }

$$
p_- < p \le  p_+.
$$
 {\it Correspondingly,   } $ q_- <  q \le q_+.  $ \par
\vspace{3mm}
{\bf Proof. A. Case $ p = p_- = 1/(1-\alpha).  $} \par

 Let us consider the function

 $$
 f_{\Delta, \theta}(x) = x^{- \Delta} \ |\log x|^{\theta}  \ I(x \in (0,1)), \ \Delta = 1 - \alpha, \ \theta = \const > -(1- \alpha).
 $$
  Then $ \left|f_{\Delta, \theta} \right|_{p_-} < \infty,  $  but it is easy to verify that $ U_{\alpha, \beta,\lambda}[f_{\Delta,\theta}] \notin L_{q_-}.  $ \par

\vspace{3mm}

{\bf Proof. B. Case $ p = p_+ = 1/(1-\alpha - \lambda).  $} \par

 Define  a function

 $$
 g(x) = x^{-( 1 - \alpha - \lambda) } \ I(x \in (1, \infty));
 $$
then

$$
\forall p > p_+ \Rightarrow g(\cdot) \in L_p,
$$
but

$$
U_{\alpha, \beta,\lambda} [g] \notin L_{q_+}, \ q_+ = q(p_+).
$$

\vspace{4mm}

 Now we formulate and prove the main result of this article.\par

\vspace{4mm}

{\bf Theorem 2.1.}

 $$
\frac{C_3(\alpha, \beta, \lambda)}{ \left[p-p_- \right]^{\kappa} } \le  K(\alpha, \beta, \lambda; p) \le
\frac{C_4(\alpha, \beta, \lambda)}{ \left[p-p_- \right]^{\kappa}}, \hspace{5mm} p \in (p_-, \ p_+],
 $$

\vspace{4mm}

 $$
  C_3(\alpha, \beta, \lambda), \  C_4(\alpha, \beta, \lambda)  \in (0,\infty).
 $$

\vspace{4mm}

{\bf Proof. Upper bound.} \par

 The upper bound follows immediately from the inequality

 $$
 K(\alpha,\beta,\lambda;p) \le \left[ \frac{\Gamma( ( 1 - 1/p - \alpha )/\kappa ) \ \Gamma( (\alpha + \beta)/\kappa)}
 {\Gamma (( 1 - 1/p + b  )/\kappa) }  \right]^{\kappa},
 $$
see \cite{Okikiolu1}, p. 213-215. \par

\vspace{4mm}

{\bf Proof. Lower bound.} \par
 Let us consider the following example.

 $$
 f_0(x) := x^{-(1 - \alpha)} \ I(x \ge 1);
 $$
then

$$
|f_0|_p \asymp c \ (p-p_-)^{ - 1/p }, \hspace{4mm} p \in (p_-, p_+].
$$
 Further, we have denoting for the values $ x > 1, \ x \to \infty $  and $ q = q(p):  $

 $$
 u_0(x) = U_{\alpha, \beta, \lambda}[f_0](x):
 $$

$$
 u_0(x) = x^{-\beta} \int_0^x \frac{y^{-\alpha} \ y^{ - 1 + \alpha }  \ dy  }{|x-y|^{\lambda}} =
 x^{ - (\beta + \lambda) } \ \int_{1/x}^1 \frac{z^{-1} \ dz }{(1-z)^{\lambda}} \sim
$$

$$
x^{ - (\beta + \lambda) } \ \int_{1/x}^1 z^{-1} \ dz =  x^{ - (\beta + \lambda) } \ \log x;
$$

$$
| u_0|_{q}^q \asymp \int_1^{\infty} x^{-q(\beta + \lambda)} \ ( \log x)^q \ dx =
\frac{\Gamma( q+1)}{ \left[ q(\beta + \lambda) - 1 \right]^{q + 1 }  };
$$

$$
|u_0|_q \asymp \frac{1}{ (q - q_-)^{ 1  + 1/q  } } \asymp   \frac{1}{ (p - p_-)^{ 1  + 1/q  } };
$$

$$
\frac{|u_0|_q}{|f|_p} \asymp \frac{C(\alpha,\beta,\lambda)}{(p - p_-)^{\kappa }}.
$$

 This completes the proof of theorem 2.1. \par

\bigskip

\section{Multidimensional case}

\vspace{3mm}

 We recall here the definition of the so-called anisotropic Lebesgue (Lebesgue-
Riesz) spaces. More detail information about this spaces see in the books of Besov O.V.,
Ilin V.P., Nikolskii S.M. \cite{Besov1}, chapter 16,17; Leoni G. \cite{Leoni1}, chapter 11; using
for us theory of operators interpolation in this spaces see in \cite{Besov1}, \cite{Bennet1}. \par

  Let $ (X_j ,A_j , μ_j, \ j = 1,2,\ldots,d) $ be measurable spaces with sigma-finite non - trivial measures
$ μ_j . $   (It is clear that in this article $ X_j = R_+  $ and $ \mu_j  $ is ordinary Lebesgue measure.)\par

 Let  $ p = \vec{p} = (p_1, p_2, ..., p_d) $ be $ d \ − \ $ dimensional vector such that $ 1 \le p_j \le \infty. $
Recall that the {\it anisotropic} Lebesgue space $ L(\vec{p}) $ consists on all the total measurable
real valued function $ f = f(x_1, x_2, . . . , x_d) = f(x) = f(\vec{x}), \ x_j \in X_j $
with finite norm $ |f|_{\vec{p} } \stackrel{def}{=} $

$$
\left( \int_{X_d} \mu_d(dx_d) \left(  \int_{X_{d-1}} \mu_{d-1}(dx_{d-1}) \ldots
\left( \int_{X_1} \mu_1(dx_1) |f(x_1,x_2, \ldots, x_d)|^{p_1}  \right)^{p_2/p_1} \right)^{p_3/p_2}  \ldots \right)^{1/p_d}.
$$

 Note that in general case $ |f|{p_1,p_2} \ne |f|{p_2,p_1}, $ but $ |f|_{p,p} = |f|_p. $
Observe also that if $ f(x_1, x_2) = g_1(x_1) · g_2(x_2), $ (condition of factorization), then
$ |f|_{p_1,p_2} = |g_1|_{p_1} · |g_2|_{p_2}, $ (formula of factorization).\par

 Let also
 $$
 \vec{\alpha} = \{ \alpha_1, \alpha_2, \ldots, \alpha_d \},  \  \vec{\beta} = \{ \beta_1, \beta_2, \ldots, \beta_d \},
 $$

$$
\vec{\lambda} = \{ \lambda_1, \lambda_2, \ldots, \lambda_d \} \eqno(3.1)
$$
be three numerical $ d \ -  $ dimensional vectors such that

$$
0 < \alpha_i, \beta_i, \lambda_i; \ \alpha_i + \beta_i + \lambda_i < 1, \ i = 1,2,\ldots,d. \eqno(3.2)
$$
 We define the multidimensional (generalized) Cesaro - Hardy operator $ U_{  \vec{\alpha},  \vec{\beta}, \vec{\lambda}}[f](\vec{x}), \
\vec{x}  \in (R_+)^d  $ as follows: $  U_{  \vec{\alpha},  \vec{\beta}, \vec{\lambda}}[f](\vec{x}) \stackrel{def}{=}  $

$$
\int_0^{x_1} \frac{x_1^{-\beta_1} \ y_1^{-\alpha_1} \ dy_1 }{|x_1 - y_1|^{\lambda_1}} \left[ \int_0^{x_2  }
\frac{x_2^{-\beta_2} \ y_2^{-\alpha_2} \ \ dy_2 }{|x_2 - y_2|^{\lambda_2}} \left[\ldots \left[ \int_0^{x_d}
\frac{x_d^{-\beta_d} \ y_d^{-\alpha_d} \ f(\vec{y}) \ d y_d }{ |x_d - y_d|^{\lambda_d}  } \right] \right] \right]. \eqno(3.3)
$$

 Let  $ p = \vec{p} = (p_1, p_2, ..., p_d) $ and $ q = \vec{q} = (q_1, q_2,\ldots, q_d)    $
 be two  $ d \ − \ $ dimensional vectors such that $ 1 < p_j, q_j < \infty. $ \par

 \vspace{4mm}

 {\it We impose on the parameters $ \{  p_j, q_j  \}  $ in this section the following condition:  }

 $$
 1 + \frac{1}{q_j} = \frac{1}{p_j} + \alpha_j + \beta_j + \lambda_j, \ j = 1,2,\ldots, d. \eqno(3.4)
 $$

  Denote

 $$
 p_-^{(j)} = \frac{1}{1 - \alpha_j}, \  p_+^{(j)} = \frac{1}{1 - \alpha_j - \lambda_j}, \eqno(3.5a)
 $$

 $$
 q_-^{(j)}  = \frac{1}{\beta_j + \lambda_j}, \  q_+^{(j)} = \frac{1}{\beta_j}, \eqno(3.5b)
 $$

 $$
 \kappa_j =  \alpha_j + \beta_j + \lambda_j. \eqno(3.6)
 $$
  The equations (3.5a) and (3.5b) uniquely define the functions $  q_j = q_j(p_j)  $  and conversely the functions
 $ p_j = p_j(q_j);  $  wherein $ p_j \in ( p_-^{(j)} p_+^{(j)} ) $  and correspondingly $ q_j \in (q_-^{(j)},  q_-^{(j)}).  $ \par

   Introduce  as before the following function:

 $$
 K_{\vec{\alpha}, \vec{\beta}, \vec{\lambda} } (\vec{p}) = \sup_{ f \in L(\vec{p}), f \ne 0 }
\left[ \frac{|U_{\vec{\alpha}, \vec{\beta}, \vec{\lambda}} [f]|_{\vec{q}}}{|f|_{\vec{p}}} \right], \
 \vec{q} = \vec{q} ( \vec{p}). \eqno(3.7)
 $$

 \vspace{4mm}

 {\bf Theorem 3.1. A.} {\it The "constant"  $  K_{\vec{\alpha}, \vec{\beta}, \vec{\lambda} } (\vec{p}) $   is finite iff  }

$$
 \forall  j = 1,2,\ldots , d \ \Rightarrow   p_-^{(j)} <  p_j  \le  p_+^{(j)} \eqno(3.8)
$$
{\it  and equation (3.4) is satisfied.}\par

\vspace{4mm}

{\bf B.} {\it If both these conditions are satisfied, then   }

$$
\frac{C_5 (\vec{\alpha}, \vec{\beta}, \vec{\lambda} )}{\prod_{j=1}^d (p_j - p_-^{(j)})^{ \kappa_j} } \le
K_{\vec{\alpha}, \vec{\beta}, \vec{\lambda} } (\vec{p}) \le
\frac{C_6 (\vec{\alpha}, \vec{\beta}, \vec{\lambda} )}{\prod_{j=1}^d (p_j - p_-^{(j)})^{ \kappa_j} }. \eqno(3.9)
$$

\vspace{4mm}

{\bf Proof} is quite similar to the analogous proof for the weight Riesz potential, see  \cite{Ostrovsky101}
and may be omitted.\par
 In particular, the example for lower estimate may be constructed as a {\it factorable  } function of a view

 $$
 f_0(\vec{x}) = \prod_{j=1}^d g_j(x_j).
 $$

\bigskip

\section{Generalization on the Grand Lebesgue Spaces (GLS). }

\vspace{3mm}

We recall first of all here  for reader conventions some definitions and facts from
the theory of GLS spaces.\par

Recently, see
\cite{Fiorenza1}, \cite{Fiorenza2},\cite{Ivaniec1}, \cite{Ivaniec2}, \cite{Jawerth1},
\cite{Karadzov1}, \cite{Kozatchenko1}, \cite{Liflyand1}, \cite{Ostrovsky1}, \cite{Ostrovsky2} etc.
 appear the so-called Grand Lebesgue Spaces GLS
 $$
 G(\psi) = G = G(\psi ; A;B);  \ A;B = \const; \ A \ge 1, \ B \le \infty
 $$
spaces consisting on all the measurable functions $ f : X \to R  $ with finite norms

$$
||f||G(\psi) \stackrel{def}{=} \sup_{p \in (A;B)} \left[\frac{|f|_p}{\psi(p)} \right]. \eqno(4.1)
$$

 Here $ \psi = \psi(p), \ p \in (A,B) $ is some continuous positive on the {\it open} interval $ (A;B) $ function such
that

$$
\inf_{p \in(A;B)} \psi(p) > 0. \eqno(4.2)
$$

We will denote
$$
\supp(\psi) \stackrel{def}{=} (A;B).
$$

The set of all such a functions with support $ \supp(\psi) = (A;B) $ will be denoted by  $  \Psi(A;B). $  \par

This spaces are rearrangement invariant; and are used, for example, in
the theory of Probability, theory of Partial Differential Equations,
 Functional Analysis, theory of Fourier series,
 Martingales, Mathematical Statistics, theory of Approximation  etc. \par

 Notice that the classical Lebesgue - Riesz spaces $ L_p $  are extremal case of Grand Lebesgue Spaces, see
 \cite{Ostrovsky2},  \cite{Ostrovsky100}. \par

 Let a function $  f:  R_+ \to R  $ be such that

 $$
 \exists (A,B): \ 1 \le A < B \le \infty \ \Rightarrow  \forall p \in (A,B) \ |f|_p < \infty.
 $$
Then the function $  \psi = \psi(p) $ may be naturally defined by the following way:

$$
\psi_f(p) := |f|_p, \ p \in (A,B). \eqno(4.3)
$$

 Let now the (measurable) function $ f: R_+ \to R, \ f \in G\psi \  $  for some $  \psi(\cdot) $ with support
 $ \supp \psi = (A,B)  $ for which

 $$
   (a,b) := (A,B) \cap (p_-, p_+) \ne \emptyset.  \eqno(4.4)
 $$

 We define a new  $ \psi \ - $ function, say $ \psi_{K} =  \psi_{K}(q)  $  as follows.

 $$
 \psi_{K}(q) =  K(\alpha, \beta, \lambda; p(q)) \cdot \psi(p(q)), \ p \in (a,b). \eqno(4.5)
 $$

\vspace{4mm}

{\bf Theorem 4.1.}     {\it Denote   }

$$
\psi_{a,b}(p) = \psi(p) \cdot  I(p \in (a,b)).
$$

{\it We assert  under condition (4.4):}

$$
||U_{\alpha, \beta, \lambda} [f] ||G \psi_K \le 1 \cdot ||f||G\psi_{a,b},  \eqno(5.6)
$$
{\it  where the constant "1" is the best possible.} \par

\vspace{4mm}

{\bf Proof. Upper bound.} \par
 Let further  in this section $ p \in (a,b). $ We can and will suppose without loss of generality  $  ||f||G\psi_{a,b} = 1. $
Then

$$
|f|_p \le \psi_{a,b}(p), \ p \in (a,b). \eqno(4.7)
$$

 We conclude after  substituting into the inequality (1.6)

 $$
 |U_{\alpha,\beta,\lambda}[f] |_{q(p)} \le K(\alpha, \beta, \lambda; p) \cdot  \psi_{a,b}(p) = \psi_K(p), \  p \in (a,b). \eqno(4.8)
 $$
The inequality (3.6) follows from  (3.8) after substitution $ p = p(q). $ \par

\vspace{4mm}

{\bf Proof. Exactness.} \par
 The exactness of the constant "1"  in the proposition (3.8) follows from  the  theorem 2.1 in the article
\cite{Ostrovsky100}. \par

\hfill $\Box$

\bigskip

\section{Concluding remarks}

\vspace{3mm}

{\bf 1.}  Analogously may be investigated the "conjugate" operator  of a view

$$
W_{\alpha, \beta, \lambda} [f](x) =
\frac{x^{ - \beta }}{\Gamma(\alpha)}
 \int_x^{\infty} \frac{ y^{-\alpha} \ f(y) dy}{  (y-x)^{\lambda}}.\eqno(5.1)
$$
see \cite{Okikiolu1}, p. 213, \cite{Mitrinovic1} , p. 173-176.\par

\vspace{3mm}

{\it We retain at the same notations and the restrictions on the parameters $ ( \alpha, \beta, \lambda, \kappa; p,q; p_-, p_+)  $
as in the second section; in opposite case }  $  K^{(W)}(\alpha,\beta,\lambda ;p) = \infty. $  \par
 In particular, again

 $$
 1 + \frac{1}{q} = \frac{1}{p} + \kappa, \ \ q =  q(p).
 $$

\vspace{3mm}

 We denote as before

 $$
 K^{(W)}(\alpha,\beta,\lambda ;p) = \sup_{f \in L(p),  f \ne 0}  \left[ \frac{|W_{\alpha,\beta,\lambda}[f] |_{q(p)}}{|f|_p}  \right]. \eqno(5.2)
 $$

 \vspace{3mm}

 {\bf Proposition 5.1.}

 $$
 K^{(W)}(\alpha, \beta, \lambda; p) \asymp \frac{C_7(\alpha, \beta, \lambda)}{(p - p_-)^{\kappa}},
  \  p \in (p_-, p_+). \eqno(5.3)
 $$

\vspace{3mm}

{\bf 2.} For the weighted convolution operator

$$
V_S[f](x) = \frac{1}{S(x)} \int_0^x s(t-x) \ f(t) \ dt,
$$
where

$$
s(t) > 0, \ S(t) = \int_0^t s(x) \ dx, \ L := \sup_{x,y: x > y} [s(x)/s(y)] < \infty,
$$
it is known \cite{Mitrinovic1}, p. 173-176 that

$$
|V_S[f]|_p \le \frac{L \ p^2}{p-1} \ |f|_p. \eqno(5.4)
$$

  Note that this operator contains as a particular case the well-known Rieman-Liouville
fractional  derivative operator, in which $  s(x) = x^{\beta - 1}, \ \beta = \const \in (0,1). $\par

\vspace{3mm}

{\bf 3.} It may be investigated analogously the discrete version of considered inequalities, i.e. when

$$
M[a](n) = \sum_{m=0}^{\infty} M(m,n) \ a(m), \ m,n = 0,1,2,\ldots \eqno(5.5)
$$
at least in the case when the function $ M(x,y), \ x,y > 0  $ is homogeneous of degree $ -1; $ see \cite{Hardy1}. \par
 Here as usually

 $$
 |\vec{a} |_p = |  \{   a(n) \} |_p = \left[ \sum_{n=0}^{\infty} |a(n)|^p  \right]^{1/p}.
 $$


\vspace{5mm}

\begin{thebibliography}{99}

\vspace{4mm}

\bibitem{Bennet1}
{\sc C. Bennet and R. Sharpley.} {\it Interpolation of operators.}
Orlando, Academic Press Inc., 1988.
\bibitem{Besov1}
{\sc O.V. Besov, V.P.Ilin, S.M.Nikolskii.}  {\it Integral representation of functions
and imbedding theorems.} Vol.2; Scripta Series in Math., V.H.Winston and Sons,
(1979), New York, Toronto, Ontario, London.
\bibitem{Carro1}
{\sc M. Carro and J. Martin. } {\it Extrapolation theory for the real
interpolation method.} Collect. Math. {\bf 33}(2002), 163--186.
\bibitem{Fiorenza1}
{\sc A. Fiorenza.} {\it Duality and reflexivity in grand Lebesgue
spaces.} Collect. Math. {\bf 51}(2000), 131--148.
\bibitem{Fiorenza2}
{\sc A. Fiorenza and G.E. Karadzhov. } {\it Grand and small Lebesgue
spaces and their analogs.} Consiglio Nationale Delle Ricerche,
Instituto per le Applicazioni del Calcoto Mauro Picine", Sezione
di Napoli, Rapporto tecnico 272/03(2005).
\bibitem{Ivaniec1}
{\sc T. Iwaniec and C. Sbordone.}  {\it On the integrability of the
 Jacobian under minimal hypotheses.} Arch. Rat.Mech. Anal., {\bf
119}(1992), 129--143.
\bibitem{Ivaniec2}
{\sc T. Iwaniec, P. Koskela and J. Onninen. } {\it Mapping of Finite
Distortion: Monotonicity and Continuity.} Invent. Math. {\bf
144}(2001), 507--531.
\bibitem{Jawerth1}
{\sc B. Jawerth and M. Milman.} {\it Extrapolation theory with
applications.} Mem. Amer. Math. Soc. {\bf 440}(1991).
\bibitem{Karadzov1}
{\sc G.E. Karadzhov and M. Milman.} {\it Extrapolation theory: new
results and applications.} J. Approx. Theory, {\bf 113}(2005),
38--99.
\bibitem{Kozatchenko1}
{\sc Yu.V. Kozatchenko and E.I. Ostrovsky.} {\it Banach spaces of random
variables of subgaussian type.} Theory Probab. Math. Stat., Kiev,
1985,  42--56 (Russian).
\bibitem{Krein1}
{\sc S.G. Krein, Yu. V. Petunin and E.M. Semenov.} {\it Interpolation of
Linear operators.} New York, AMS, 1982.
\bibitem{Liflyand1}
{\sc E. Liflyand, E. Ostrovsky, L. Sirota.} Structural Properties of Bilateral Grand Lebesgue Spaces.
Turk. J. Math.; {\bf 34} (2010), 207-219.
\bibitem{Leoni1}
{\sc G.Leoni.}   {\it A first Course in Sobolev Spaces.} Graduate Studies in Mathematics,
v. 105, AMS, Providence, Rhode Island, (2009).
\bibitem{Mitrinovic1}
{\sc D.S. Mitrinovic, J.E.Pecaric, and A.M.Fink.} {\it Inequalities  Involving Functions and Their Integrals
and Derivatives. } Kluvner Verlag, (1991).
\bibitem{Ostrovsky1}
{\sc E.I. Ostrovsky.} {\it Exponential Estimations for Random Fields.}
Moscow - Obninsk, OINPE, 1999 (Russian).
\bibitem{Ostrovsky2}
{\sc E. Ostrovsky and L.Sirota.} {\it Moment Banach spaces: theory and applications.}
HAIT Journal of Science and Engeneering, {\bf C}, Volume 4, Issues 1-2,
pp. 233 - 262, (2007).
\bibitem{Ostrovsky100}
{\sc E. Ostrovsky and L.Sirota.} {\it  Boundedness of operators in bilateral Grand Lebesgue Spaces,
with exact and weakly exact constant calculation. }
arXiv:1104.2963v1 [math.FA] 15 Apr 2011
\bibitem{Ostrovsky101}
{\sc E. Ostrovsky and L.Sirota.}
{\it  Multiple weight Riesz and Fourier transforms  in bilateral anosotropic Grand
Lebesgue Spaces.  }
arXiv:1208.2392v1 [math.FA] 12 Aug 2012
\bibitem{Perez1}
{\sc C.Perez (joint work  with A.Lerner, S.Ombrosi, K.Moen and  R.Torres).} {\it   Sharp
Weighted Bound for Zygmund Singular Integral Operators and Sobolev Inequalities.}
In: "Oberwolfach Reports", Vol. Nunber 3, p. 1828 - 1830; EMS Publishing House,
ETH - Zentrum FLIC1, CH - 8092, Zurich, Switzerland.
\bibitem{Stein1}
{\sc E.M.Stein.} {\it Singular Integrals and Differentiability Properties  of Functions. }
Princeton University Press, Princeton, (1992).

\vspace{5mm}

\bibitem{Fujita1}
{\sc H.Fujita and T. Kato. } {\it On the Navier-Stokes initial value problem I. } Arch. Ration. Mech.
Anal., 16(1964), 269 \ – \ 315.
\bibitem{Giga4}
{\sc Y.Giga and  H.Sohr. }  {\it Abstract $ L^p \ - $ estimates for the Cauchy  problem with Applications to the
Navier-Stokes equations in  exteroir domains. } J. Func. Anal., {\bf 102} (1991), 72 \ -  \ 94.
\bibitem{Hardy1}
G.H.Hardy, J.E. Litlewood and G.P\'olya. {\it Inequalities.}
Cambridge, University Press (1952).
\bibitem{Kato1}
{\sc T. Kato. }  {\it Strong $ L_p $ solutions of the Navier-Stokes equations in $ R^m $ with applications to
weak solutions. } Math. Z., 187 \ (1984), 471 \ – \ 480.
\bibitem{Lizama1}
{\sc C.Lizama, P.J.Miana, R.Ponce and L.S.Anchez-Lajusticia. } {\it  On the boundedness of generalized Cesaro
 operators on Sobolev spaces.}
arXiv:1304.1622v1 [math.FA] 5 Apr 2013
\bibitem{Okikiolu1}
{\sc G.O.Okikiolu.} {\it Aspects of the theory of bounded Integral Operators in the $ L^p $  Spaces. }
Academic Press; London,   New Yotk; (1971).
\bibitem{Ostrovsky3}
{\sc E.Ostrovsky,  L.Sirota, E.Rogover.} {\it  Riesz's and Bessel's operators in Bilateral Grand Lebesgue Spaces. }
arXiv:0907.3321v1 [math.FA] 19 Jul 2009
\bibitem{Ostrovsky4}
{\sc E.Ostrovsky,  L.Sirota.} {\it Weight Hardy-Littlewood inequalities for different powers. }
arXiv:0910.5880v1 [math.FA] 30 Oct 2009
\bibitem{Ostrovsky5}
{\sc E.Ostrovsky,  L.Sirota.}
{\it  Quantitative lower bounds for lifespan for solution of Navier-Stokes equations.  }
arXiv:1306.6211v1 [math.AP] 26 Jun 2013
\bibitem{Talenti1}
{\sc G.Talenti.} {\it Inequalities in Rearrangement Invariant Function Spaces. Nonlinear Analysis, Function
 Spaces and Applications.} Prometheus, Prague, {\bf 5}, (1995), 177-230.

\end{thebibliography}
\end{document}